\documentclass[12pt]{article}

\usepackage{amsfonts,amssymb,amsthm,amsmath,latexsym}
\usepackage{subeqnarray,eqsection,indent}
\usepackage{cite}

\input cyracc.def
\font\tencyr=wncyr10
\def\cyr{\tencyr\cyracc}

\begin{document}

\title{A Ridiculously Simple and Explicit \\ Implicit Function Theorem}

\author{
     {\small Alan D.~Sokal\thanks{Also at Department of Mathematics,
           University College London, London WC1E 6BT, England.}}  \\[-2mm]
     {\small\it Department of Physics}       \\[-2mm]
     {\small\it New York University}         \\[-2mm]
     {\small\it 4 Washington Place}          \\[-2mm]
     {\small\it New York, NY 10003 USA}      \\[-2mm]
     {\small\tt sokal@nyu.edu}               \\[-2mm]
     {\protect\makebox[5in]{\quad}}  
     \\
}

\date{January 31, 2009 \\[5mm]
      \emph{Dedicated to the memory of Pierre Leroux}}

\maketitle
\thispagestyle{empty}   
\begin{abstract}
I show that the general implicit-function problem
(or parametrized fixed-point problem)
in one complex variable has an explicit series solution
given by a trivial generalization of the Lagrange inversion formula.
I give versions of this formula for both
analytic functions and formal power series.
\end{abstract}


\bigskip
\noindent
{\bf Mathematics Subject Classification (MSC 2000) codes:}
30B10 (Primary); 05A15, 13F25, 26B10, 30A20, 32A05, 47J07 (Secondary).

\clearpage

\newtheorem{theorem}{Theorem}[section]
\newtheorem{proposition}[theorem]{Proposition}
\newtheorem{lemma}[theorem]{Lemma}
\newtheorem{corollary}[theorem]{Corollary}
\newtheorem{definition}[theorem]{Definition}
\newtheorem{conjecture}[theorem]{Conjecture}
\newtheorem{question}[theorem]{Question}
\newtheorem{example}[theorem]{Example}
\newtheorem{remark}[theorem]{Remark}

\renewcommand{\theenumi}{\alph{enumi}}
\renewcommand{\labelenumi}{(\theenumi)}
\def\eop{\hbox{\kern1pt\vrule height6pt width4pt
depth1pt\kern1pt}\medskip}
\def\prf{\par\noindent{\bf Proof.\enspace}\rm}
\def\rmk{\par\medskip\noindent{\bf Remark\enspace}\rm}

\newcommand{\be}{\begin{equation}}
\newcommand{\ee}{\end{equation}}
\newcommand{\<}{\langle}
\renewcommand{\>}{\rangle}
\newcommand{\widebar}{\overline}
\def\reff#1{(\protect\ref{#1})}
\def\spose#1{\hbox to 0pt{#1\hss}}
\def\ltapprox{\mathrel{\spose{\lower 3pt\hbox{$\mathchar"218$}}
    \raise 2.0pt\hbox{$\mathchar"13C$}}}
\def\gtapprox{\mathrel{\spose{\lower 3pt\hbox{$\mathchar"218$}}
    \raise 2.0pt\hbox{$\mathchar"13E$}}}
\def\textprime{${}^\prime$}
\def\proof{\par\medskip\noindent{\sc Proof.\ }}
\def\firstproof{\par\medskip\noindent{\sc First Proof.\ }}
\def\secondproof{\par\medskip\noindent{\sc Second Proof.\ }}
\def\thirdproof{\par\medskip\noindent{\sc Third Proof.\ }}
\def\qed{ $\square$ \bigskip}
\def\proofof#1{\bigskip\noindent{\sc Proof of #1.\ }}
\def\firstproofof#1{\bigskip\noindent{\sc First Proof of #1.\ }}
\def\secondproofof#1{\bigskip\noindent{\sc Second Proof of #1.\ }}
\def\thirdproofof#1{\bigskip\noindent{\sc Third Proof of #1.\ }}
\def\half{ {1 \over 2} }
\def\third{ {1 \over 3} }
\def\twothird{ {2 \over 3} }
\def\smfrac#1#2{{\textstyle{#1\over #2}}}
\def\smhalf{ {\smfrac{1}{2}} }
\newcommand{\real}{\mathop{\rm Re}\nolimits}
\renewcommand{\Re}{\mathop{\rm Re}\nolimits}
\newcommand{\imag}{\mathop{\rm Im}\nolimits}
\renewcommand{\Im}{\mathop{\rm Im}\nolimits}
\newcommand{\sgn}{\mathop{\rm sgn}\nolimits}
\newcommand{\tr}{\mathop{\rm tr}\nolimits}
\newcommand{\supp}{\mathop{\rm supp}\nolimits}
\def\hboxscript#1{ {\hbox{\scriptsize\em #1}} }
\renewcommand{\emptyset}{\varnothing}

\newcommand{\restrict}{\upharpoonright}

\newcommand{\scra}{{\mathcal{A}}}
\newcommand{\scrb}{{\mathcal{B}}}
\newcommand{\scrc}{{\mathcal{C}}}
\newcommand{\scre}{{\mathcal{E}}}
\newcommand{\scrf}{{\mathcal{F}}}
\newcommand{\scrg}{{\mathcal{G}}}
\newcommand{\scrh}{{\mathcal{H}}}
\newcommand{\scrk}{{\mathcal{K}}}
\newcommand{\scrl}{{\mathcal{L}}}
\newcommand{\scro}{{\mathcal{O}}}
\newcommand{\scrp}{{\mathcal{P}}}
\newcommand{\scrr}{{\mathcal{R}}}
\newcommand{\scrs}{{\mathcal{S}}}
\newcommand{\scrt}{{\mathcal{T}}}
\newcommand{\scrv}{{\mathcal{V}}}
\newcommand{\scrw}{{\mathcal{W}}}
\newcommand{\scrz}{{\mathcal{Z}}}

\newcommand{\ahat}{{\widehat{a}}}
\newcommand{\Zhat}{{\widehat{Z}}}
\renewcommand{\k}{{\mathbf{k}}}
\newcommand{\n}{{\mathbf{n}}}
\newcommand{\vv}{{\mathbf{v}}}
\newcommand{\bv}{{\mathbf{v}}}
\newcommand{\w}{{\mathbf{w}}}
\newcommand{\x}{{\mathbf{x}}}
\newcommand{\g}{{\boldsymbol{g}}}
\newcommand{\cc}{{\mathbf{c}}}
\newcommand{\zero}{{\mathbf{0}}}
\newcommand{\one}{{\mathbf{1}}}
\newcommand{\balpha}{{\boldsymbol{\alpha}}}

\newcommand{\C}{{\mathbb C}}
\newcommand{\D}{{\mathbb D}}
\renewcommand{\Dbar}{{\overline{\mathbb D}}}
\newcommand{\Z}{{\mathbb Z}}
\newcommand{\N}{{\mathbb N}}
\newcommand{\Q}{{\mathbb Q}}
\newcommand{\R}{{\mathbb R}}
\newcommand{\RR}{{\mathbb R}}

\def\cbar{{\overline{C}}}
\def\ctilde{{\widetilde{C}}}
\def\zbar{{\overline{Z}}}
\def\pitilde{{\widetilde{\pi}}}

%
%
\newcommand{\stirlingsubset}[2]{\genfrac{\{}{\}}{0pt}{}{#1}{#2}}
\newcommand{\stirlingcycle}[2]{\genfrac{[}{]}{0pt}{}{#1}{#2}}
\newcommand{\assocstirlingsubset}[3]{{\genfrac{\{}{\}}{0pt}{}{#1}{#2}}_{\! \ge #3}}
\newcommand{\assocstirlingcycle}[3]{{\genfrac{[}{]}{0pt}{}{#1}{#2}}_{\ge #3}}
\newcommand{\genstirlingsubset}[4]{{\genfrac{\{}{\}}{0pt}{}{#1}{#2}}_{\! #3,#4}}
\newcommand{\euler}[2]{\genfrac{\langle}{\rangle}{0pt}{}{#1}{#2}}
\newcommand{\eulergen}[3]{{\genfrac{\langle}{\rangle}{0pt}{}{#1}{#2}}_{\! #3}}
\newcommand{\eulersecond}[2]{\left\langle\!\! \euler{#1}{#2} \!\!\right\rangle}
\newcommand{\eulersecondgen}[3]{{\left\langle\!\! \euler{#1}{#2} \!\!\right\rangle}_{\! #3}}
\newcommand{\binomvert}[2]{\genfrac{\vert}{\vert}{0pt}{}{#1}{#2}}


\newenvironment{sarray}{
             \textfont0=\scriptfont0
             \scriptfont0=\scriptscriptfont0
             \textfont1=\scriptfont1
             \scriptfont1=\scriptscriptfont1
             \textfont2=\scriptfont2
             \scriptfont2=\scriptscriptfont2
             \textfont3=\scriptfont3
             \scriptfont3=\scriptscriptfont3
           \renewcommand{\arraystretch}{0.7}
           \begin{array}{l}}{\end{array}}

\newenvironment{scarray}{
             \textfont0=\scriptfont0
             \scriptfont0=\scriptscriptfont0
             \textfont1=\scriptfont1
             \scriptfont1=\scriptscriptfont1
             \textfont2=\scriptfont2
             \scriptfont2=\scriptscriptfont2
             \textfont3=\scriptfont3
             \scriptfont3=\scriptscriptfont3
           \renewcommand{\arraystretch}{0.7}
           \begin{array}{c}}{\end{array}}

\clearpage

\section{Introduction}

It is well known to both complex analysts and combinatorialists
that the problem of reverting a power series,
i.e.\ solving $f(z) = w$, has an explicit solution known as the
Lagrange (or Lagrange--B\"urmann) inversion formula
\cite{Whittaker_27,Caratheodory_58,Markushevich_65,Henrici_74,%
Goulden_83,Wilf_94,Bergeron_98,Stanley_99}.
What seems to be less well known is that the more general
implicit-function problem $F(z,w) = 0$ also has a simple explicit solution,
given by Yuzhakov \cite{Yuzhakov_75a} in 1975
(see also \cite{Yuzhakov_75b,Bolotov_78,Aizenberg_83}).
My purpose here is to give a slightly more general and flexible version
of Yuzhakov's formula, and to show that its proof is an utterly trivial
generalization of the standard proof of the Lagrange inversion formula.
A formal-power-series version of the formula presented here appears also in
\cite[Exercise~5.59, pp.~99 and 148]{Stanley_99},
but its importance for the implicit-function problem
does not seem to be sufficiently stressed.

Let us begin by recalling the Lagrange inversion formula:
if $f(z) = \sum_{n=1}^\infty a_n z^n$ with $a_1 \neq 0$
(interpreted either as an analytic function or as a formal power series),
then
\begin{equation}
   f^{-1}(w)  \;=\;
   \sum_{m=1}^\infty {w^m \over m} \,
      [\zeta^{m-1}] \biggl( {\zeta \over f(\zeta)}\biggr)^m
   \;,
\end{equation}
where $[\zeta^n] g(\zeta) = g^{(n)}(0)/n!$ denotes the coefficient of $\zeta^n$
in the power series $g(\zeta)$.
More generally, if $h(z) = \sum_{n=0}^\infty b_n z^n$, we have
\begin{equation}
   h(f^{-1}(w))  \;=\;
   h(0) \,+\,
   \sum_{m=1}^\infty {w^m \over m} \,
      [\zeta^{m-1}] \, h'(\zeta) \biggl( {\zeta \over f(\zeta)}\biggr)^m
   \;.
\end{equation}
Proofs of these formulae can be found in many books on analytic function
theory \cite{Whittaker_27,Caratheodory_58,Markushevich_65,Henrici_74}
and enumerative combinatorics \cite{Goulden_83,Wilf_94,Bergeron_98,Stanley_99}.

It is convenient to introduce the function (or formal power series)
$g(z) = z/f(z)$;  then the equation $f(z)=w$ can be rewritten as $z = g(z) w$,
and its solution $z = \varphi(w) = f^{-1}(w)$ is given by the power series
\begin{equation}
   \varphi(w)  \;=\;
   \sum_{m=1}^\infty {w^m \over m} \, [\zeta^{m-1}] g(\zeta)^m
 \label{eq.star1}
\end{equation}
and
\begin{equation}
   h(\varphi(w))  \;=\;
   h(0) \,+\
   \sum_{m=1}^\infty {w^m \over m} \, [\zeta^{m-1}] h'(\zeta) g(\zeta)^m
   \;.
 \label{eq.star2}
\end{equation}
There is also an alternate form
\begin{equation}
   h(\varphi(w))  \;=\;
   h(0) \,+\
   \sum_{m=1}^\infty w^m  \,
      [\zeta^m] h(\zeta) \bigl[ g(\zeta)^m \,-\, \zeta g'(\zeta) g(\zeta)^{m-1}
                         \bigr]
   \;.
 \label{eq.star2a}
\end{equation}

Consider now the more general problem of solving $z = G(z,w)$,
where $G(0,0) = 0$ and $|(\partial G/\partial z)(0,0)| < 1$.
I shall prove here that its solution $z = \varphi(w)$ is given by
the function series
\begin{equation}
   \varphi(w)  \;=\;
   \sum_{m=1}^\infty {1 \over m} \, [\zeta^{m-1}] G(\zeta,w)^m
   \;.
 \label{eq.star3}
\end{equation}
More generally, for any analytic function (or formal power series)
$H(z,w)$, we have
\begin{eqnarray}
   H(\varphi(w),w)
   & = &
   H(0,w) \,+\,
   \sum_{m=1}^\infty {1 \over m} \,
      [\zeta^{m-1}] {\partial H(\zeta,w) \over \partial\zeta} \, G(\zeta,w)^m
            \label{eq.star4} \\[2mm]
   & = &
   H(0,w) \,+\,
   \sum_{m=1}^\infty 
      [\zeta^m] H(\zeta,w)  \Bigl[ G(\zeta,w)^m \:-\:
                     \zeta \, {\partial G(\zeta,w) \over \partial\zeta} \,
                     G(\zeta,w)^{m-1}
                              \Bigr]
   \;.
            \nonumber \\ \label{eq.star4a}
\end{eqnarray}
The formulae \reff{eq.star3}--\reff{eq.star4a} are manifestly generalizations
of the Lagrange inversion formulae \reff{eq.star1}--\reff{eq.star2a},
to which they reduce when $G(z,w) = g(z) w$ and $H(z,w) = h(z)$.
It turns out that the proof of \reff{eq.star3}--\reff{eq.star4a} is, likewise,
virtually a verbatim copy of the standard proof of
\reff{eq.star1}--\reff{eq.star2a}:
the variables $w$ simply ``go for the ride''.

The problem of solving $z = G(z,w)$ can alternatively be interpreted
as a fixed-point problem for the family of maps $z \mapsto G(z,w)$
parametrized by $w$.
{}From this point of view, \reff{eq.star3}--\reff{eq.star4a} are simply
a function series giving the unique solution of this fixed-point problem
under a suitable ``Rouch\'e-contraction'' hypothesis
[see hypothesis (c) of Theorem~\ref{thm1} below].
Once again, the variables $w$ simply ``go for the ride''.

Before proving \reff{eq.star3}--\reff{eq.star4a},
let us observe how these formulae solve the implicit-function problem
$F(z,w) = 0$, where $F(0,0) = 0$ and
$(\partial F/\partial z)(0,0) \equiv a_{10} \neq 0$.
%
%
It suffices to choose {\em any}\/ analytic function $\gamma(z,w)$
satisfying $\gamma(0,0) \neq 0$, and then define
\begin{equation}
   G(z,w)  \;=\;  z \,-\,  \gamma(z,w) F(z,w)
   \;.
 \label{eq.star6}
\end{equation}
Clearly $F(z,w) = 0$ is equivalent to $z = G(z,w)$,
at least locally in a neighborhood of $(z,w) = (0,0)$.
Then $(\partial G/\partial z)(0,0) = 1 - \gamma(0,0) a_{10}$;
so for \reff{eq.star3}--\reff{eq.star4a} to be applicable,
it suffices to arrange that $|1 - \gamma(0,0) a_{10}| < 1$,
which can easily be done by a suitable choice of $\gamma(0,0)$
[namely, by choosing $\gamma(0,0)$ to lie in the open disc
 of radius $1/|a_{10}|$ centered at $1/a_{10}$].
I wish to stress that each such choice of a function $\gamma$ gives rise to
a valid but different expansion \reff{eq.star3}--\reff{eq.star4a}
for the solution $z = \varphi(w)$.
Even in the special case of the inverse-function problem $f(z) = w$,
this flexibility exists and is useful (see Example~\ref{exam3.2} below).
One important class of choices has $\gamma(0,0) = 1/a_{10}$,
so that $(\partial G/\partial z)(0,0) = 0$;
this latter condition leads to a slight simplification in the formulae
(see Remark~\ref{remark.calpha} below).
This special class in turn contains two important subclasses:
\begin{itemize}
  \item Yuzhakov \cite{Yuzhakov_75a} takes $\gamma$ to be
     the constant function $1/a_{10}$.\footnote{
           His subsequent generalization \cite{Yuzhakov_75b}
           \cite[Theorem~20.2 and Proposition~20.4]{Aizenberg_83}
           in effect allows a preliminary change of variables
           $z' = z - \psi(w)$ for arbitrary $\psi$ satisfying $\psi(0) = 0$.
}

  \item Alternatively, we can choose $\gamma$ so that $\gamma(z,0) = z/F(z,0)$
     (this definition can still be extended to $w \neq 0$ in
      many different ways).
     Then $G(z,0) \equiv 0$, so that $G(z,w)$ has an overall factor $w$.
     This leads to a further slight simplification 
     (see again Remark~\ref{remark.calpha}).
\end{itemize}

Conversely, the problem of solving $z = G(z,w)$ is of course equivalent
to the problem of solving $\widetilde{F}(z,w) = 0$ if we set
\begin{equation}
   \widetilde{F}(z,w)  \;=\;  z \,-\,  G(z,w)
   \;;
\end{equation}
and the condition $(\partial \widetilde{F}/\partial z)(0,0) \neq 0$
is satisfied whenever $(\partial G/\partial z)(0,0) \neq 1$.
So our parametrized fixed-point problem has exactly
the same level of generality as the implicit-function problem.

The plan of this paper is as follows:
First I shall state and prove four versions of the formulae
\reff{eq.star3}--\reff{eq.star4a}:
one in terms of analytic functions (Theorem~\ref{thm1}),
and three in terms of formal power series
(Theorems~\ref{thm2.first}, \ref{thm2} and \ref{thm2bis}).
Then I shall give some examples and make some final remarks.

Since this paper is aimed at a diverse audience of analysts
and combinatorialists, I have endeavored to give more detailed proofs
than would otherwise be customary.  I apologize in advance to experts
for occasionally boring them with elementary observations.

\section{Implicit function formula: Analytic version}

In the analytic version of \reff{eq.star3}--\reff{eq.star4a},
the variable $w$ simply ``goes for the ride'';
consequently, $w$ need not be assumed to lie in $\C$,
but can lie in the multidimensional complex space $\C^M$
or even in a general topological space $W$.
I shall begin with a simple auxiliary result (Proposition~\ref{prop1})
that clarifies the meaning of the hypotheses of Theorem~\ref{thm1}.
I denote open and closed discs in $\C$ by the notations
$\D_R = \{ z \in \C \colon\, |z| < R \}$
and $\Dbar_R = \{ z \in \C \colon\, |z| \le R \}$.

\begin{proposition}
 \label{prop1}
Let $W$ be a topological space, let $R > 0$,
and let $G \colon\, \D_R \times W \to \C$
be a jointly continuous function with the property that
$G( \,\cdot\,, w)$ is analytic on $\D_R$ for each fixed $w \in W$.
Suppose further that for some $w_0 \in W$ we have
$G(0,w_0) = 0$ and $|(\partial G/\partial z)(0,w_0)| < 1$.
Then
for all sufficiently small $\rho > 0$ and $\epsilon > 0$
there exists an open neighborhood $V_\rho \ni w_0$ such that
$|G(z,w)| \le (1-\epsilon) |z|$ whenever $|z|=\rho$ and $w \in V_\rho$.
\end{proposition}

\proof
The function
\begin{equation}
   f(z)  \;=\;
   \begin{cases}
      G(z,w_0)/z                      & \text{if $0 < |z| < R$} \\[0.5mm]
      (\partial G/\partial z)(0,w_0)  & \text{if $z=0$}
   \end{cases}
\end{equation}
is analytic on $\D_R$, hence continuous on $\D_R$;
and by hypothesis $|f(0)| < 1$.
It follows that for all sufficiently small $\rho > 0$ and $\epsilon > 0$,
we have $|f(z)| \le 1-2\epsilon$ whenever $|z| \le \rho$.
We now use the following simple topological fact:
\begin{lemma}
   \label{lemma.topological}
If $F \colon X \times Y \to \R$ is continuous and $X$ is compact,
then the function $g \colon Y \to \R$ defined by
$g(y) = \sup\limits_{x \in X} F(x,y)$ is continuous.
\end{lemma}
\noindent
Applying this with $X = \{ z \in \C \colon\, |z| = \rho\}$ and $Y=V$,
we conclude that there exists an open neighborhood $V_\rho \ni w_0$ such that
$|G(z,w)/z| \le 1-\epsilon$ whenever $|z| = \rho$ and $w \in V_\rho$.
\qed

For completeness, let us give a proof of the lemma:

\proofof{Lemma~\ref{lemma.topological}}
First of all, the compactness of $X$ guarantees that
$g$ is everywhere finite.
Now, since the supremum of any family of continuous functions
is lower semicontinuous, it suffices to prove that $g$ is
upper semicontinuous, i.e.\ that for any $y_0 \in Y$
and any $\epsilon > 0$ there exists an open neighborhood $V \ni y_0$
such that $g(y) < g(y_0) + \epsilon$ for all $y \in V$.
To see this, first choose, for each $x \in X$,
open neighborhoods $U_x \ni x$ and $V_x \ni y_0$
such that $F(x',y) < F(x_,y_0) + \epsilon \le g(y_0) + \epsilon$
whenever $x' \in U_x$ and $y \in V_x$.
By compactness, there exists a finite set $\{x_1,\ldots,x_n\}$
such that $\{U_{x_i}\}_{i=1}^n$ covers $X$.
Setting $V = \bigcap\limits_{i=1}^n V_{x_i}$ gives the required neighborhood.
\qed

It is easy to see that the conclusion of Lemma~\ref{lemma.topological}
need not hold if $X$ is noncompact.
For instance, with $X = \R$ and $Y = [0,1]$,
take $F(x,y) = xy$ or $F(x,y) = \tanh(xy)$.

We can now state the principal result of this section:

\begin{theorem}[Implicit function formula --- analytic version]
 \label{thm1}
\quad\par\noindent
Let $V$ be a topological space, let $\rho > 0$,
and let $G,H \colon\; \Dbar_\rho \times V \to \C$ be functions satisfying
\begin{itemize}
   \item[(a)] $G$, $\partial G/\partial z$ and $H$ are jointly continuous on
        $\Dbar_\rho \times V$;
   \item[(b)] $G( \,\cdot\,, w)$ and $H( \,\cdot\,, w)$
        are analytic on $\D_\rho$ for each fixed $w \in V$; and
   \item[(c)] $|G(z,w)| < |z|$ whenever $|z|=\rho$ and $w \in V$.
\end{itemize}
Then for each $w \in V$, there exists a unique $z \in \D_\rho$
satisfying $z = G(z,w)$.  Furthermore, this $z = \varphi(w)$
depends continuously on $w$ and is given explicitly by the function series
\begin{eqnarray}
   H(\varphi(w),w)
   & = &
   H(0,w) \,+\,
   \sum_{m=1}^\infty {1 \over m} \,
      [\zeta^{m-1}] {\partial H(\zeta,w) \over \partial\zeta} \, G(\zeta,w)^m
            \label{eq.star4bis} \\[2mm]
   & = &
   H(0,w) \,+\,
   \sum_{m=1}^\infty
      [\zeta^m] H(\zeta,w)  \Bigl[ G(\zeta,w)^m \:-\:
                     \zeta \, {\partial G(\zeta,w) \over \partial\zeta} \,
                     G(\zeta,w)^{m-1}
                              \Bigr]
            \nonumber \\ \label{eq.star4abis}
\end{eqnarray}
which are absolutely convergent on $V$, locally uniformly on $V$.\footnote{
   ``Locally uniformly'' means that for each $w \in V$
   there exists an open neighborhood $U \ni w$
   on which the convergence is uniform.
   This implies, in particular, that the convergence is uniform
   on compact subsets of $V$
   (and is equivalent to it if $V$ is locally compact).
}

If, in addition, $V$ is an open subset of $\C^M$
and $G$ is analytic on $\D_\rho \times V$,
then $\varphi$ is analytic on $V$.

Finally, if $V$ is a polydisc
(or more generally a complete Reinhardt domain)
centered at $0 \in \C^M$,
and $G$ and $H$ are analytic on $\D_\rho \times V$
and satisfy $G(0,0) = 0$,
then $\varphi(w)$ is also given by the Taylor series
\begin{equation}
   H(\varphi(w),w)  \;=\;  \sum_{\balpha \in \N^M} c_\balpha \, w^\balpha
\end{equation}
which is absolutely convergent on $V$, uniformly on compact subsets of $V$;
here the coefficients $c_\balpha$ are given by the absolutely convergent sums
\begin{eqnarray}
   c_\balpha
   & = &
   [w^\balpha] H(0,w) \,+\,
   \sum_{m=1}^\infty {1 \over m} \,
      [\zeta^{m-1} w^\balpha]
      {\partial H(\zeta,w) \over \partial\zeta} \, G(\zeta,w)^m
            \label{eq.calpha} \\[2mm]
   & = &
   [w^\balpha] H(0,w) \,+\,
   \sum_{m=1}^\infty
      [\zeta^m w^\balpha] H(\zeta,w)  \Bigl[ G(\zeta,w)^m \:-\:
                     \zeta \, {\partial G(\zeta,w) \over \partial\zeta} \,
                     G(\zeta,w)^{m-1}
                              \Bigr]
   \;.
            \nonumber \\ \label{eq.calpha.a}
\end{eqnarray}
\end{theorem}

\medskip

\begin{remark}
  \label{remark.calpha}
\rm
In many cases \reff{eq.calpha} is actually a finite sum.
For instance, if $(\partial G/\partial z)(0,0) = 0$ ---
as occurs in particular in
Yuzhakov's \cite{Yuzhakov_75a,Yuzhakov_75b,Aizenberg_83} approach ---
then each factor of $G(\zeta,w)$ brings either $\zeta^2$ or $w$ (at least);
simple algebra then shows that the summand in \reff{eq.calpha}
is nonvanishing only for $m \le 2 |\balpha| - 1$,
where $|\balpha| = \sum_{i=1}^M \alpha_i$.\footnote{
   Yuzhakov \cite[Proposition~1]{Yuzhakov_75a}
            \cite[Proposition~20.4]{Aizenberg_83}
   writes $m \le 2 |\balpha|$,
   which is correct but not as strong as it should be.
}
Under the stronger hypothesis $G(z,0) \equiv 0$,
each factor of $G(\zeta,w)$ brings at least one $w$,
so the summand is nonvanishing only for $m \le |\balpha|$.

Analogous comments hold for \reff{eq.calpha.a},
where the conditions are $m \le 2 |\balpha|$ and $m \le |\balpha|$,
respectively.
\qed
\end{remark}

As previously stated, the proof of Theorem~\ref{thm1}
is a trivial modification of the standard textbook proof
of the Lagrange inversion formula
\cite{Caratheodory_58,Markushevich_65,Whittaker_27},
but for completeness let us give it in detail.

\proof
Hypotheses (a)--(c) combined with Rouch\'e's theorem
imply that for each $w \in V$, the number of roots (including multiplicity)
of $z - G(z,w) = 0$ in the disc $|z| < \rho$
is the same as the number of roots of $z=0$ in this disc, namely one;
so let us call this unique (and simple) root $z = \varphi(w)$.\footnote{
   It is well known (and is easily proved using Rouch\'e's theorem
   and Lemma~\ref{lemma.topological})
   that this root depends continuously on $w$.
   This will also follow from the explicit formula
   \reff{eq.star4bis}/\reff{eq.star4abis}.
}
It follows that for each $w \in V$, the function
\begin{equation}
   H(\zeta,w) \, {1 - {\partial G \over \partial\zeta}(\zeta,w)
                  \over
                  \zeta - G(\zeta,w)
                 }
\end{equation}
is continuous on $|\zeta| \le \rho$ and analytic in $|\zeta| < \rho$
except for a simple pole at $\zeta = \varphi(w)$
with residue $H(\varphi(w),w)$.\footnote{
   It is in this step that we use the continuity of
   $\partial G/\partial z$ (as well as that of $G$ and $H$)
   on the closed disc $\Dbar_\rho$.
}
Cauchy's integral formula therefore gives
\begin{equation}
   H(\varphi(w),w)
   \;=\;
   {1 \over 2\pi i} \oint\limits_{|\zeta| = \rho}
   H(\zeta,w) \, {1 - {\partial G \over \partial\zeta}(\zeta,w)
                  \over
                  \zeta - G(\zeta,w)
                 }
   \; d\zeta  \;.
 \label{eq.cauchy}
\end{equation}
Let us now expand
\begin{equation}
   {1 \over \zeta - G(\zeta,w)}
   \;=\;
   \sum_{m=0}^\infty  {G(\zeta,w)^m \over \zeta^{m+1}}
 \label{eq.neumann}
\end{equation}
and observe that this sum is absolutely convergent,
uniformly for $\zeta$ on the circle $|\zeta| = \rho$ of integration
and locally uniformly for $w \in V$
[this follows from hypothesis (c) and Lemma~\ref{lemma.topological}].
We therefore have
\begin{equation}
   H(\varphi(w),w)
   \;=\;
   \sum_{m=0}^\infty
   {1 \over 2\pi i} \oint\limits_{|\zeta| = \rho}
   {H(\zeta,w) \over \zeta^{m+1}} \,
   \biggl[ 1 - {\partial G \over \partial\zeta}(\zeta,w) \biggr] \,
   G(\zeta,w)^m \, d\zeta  \;.
 \label{eq.contour}
\end{equation}
By the Cauchy integral formula, this gives
\begin{eqnarray}
   \!\!\!
   H(\varphi(w),w)
   & = &
   \sum_{m=0}^\infty [\zeta^m] H(\zeta,w) G(\zeta,w)^m
   \:-\:
   \sum_{m=0}^\infty [\zeta^m] H(\zeta,w) \, {\partial G \over \partial\zeta} \,
                  G(\zeta,w)^m
   \qquad\qquad
        \nonumber \\
   & = &
   H(0,w) \,+\,
   \sum_{m=1}^\infty [\zeta^m] H(\zeta,w) G(\zeta,w)^m
   \,-\,
   \sum_{m=1}^\infty [\zeta^{m-1}] H(\zeta,w) \,
                                   {\partial G \over \partial\zeta} \,
                                   G(\zeta,w)^{m-1}
     \;.
     \hspace*{-1in}
     \nonumber \\
\end{eqnarray}
This proves the alternate formula \reff{eq.star4a}/\reff{eq.star4abis},
in which the sum is absolutely convergent on $V$,
locally uniformly on $V$.

To prove \reff{eq.star4}/\reff{eq.star4bis},
we start from \reff{eq.contour} and prepare an integration by parts:
\begin{eqnarray}
   {H(\zeta,w) \over \zeta^{m+1}} \, {\partial G \over \partial \zeta} \, G^m
   & = &
   {1 \over m+1} \, {H(\zeta,w) \over \zeta^{m+1}} \,
        {\partial \over \partial \zeta} (G^{m+1})
   \nonumber \\[2mm]
   & = &
   {\partial \over \partial \zeta}
     \biggl( {1 \over m+1} \, {H(\zeta,w) \over \zeta^{m+1}} \, G^{m+1} \biggr)
   \,-\,
   {G^{m+1} \over m+1} \, {\partial \over \partial \zeta}
     \biggl( {H(\zeta,w) \over \zeta^{m+1}} \biggr)
   \nonumber \\[2mm]
   & = &
   {\partial \over \partial \zeta}
     \biggl( {1 \over m+1} \, {H(\zeta,w) \over \zeta^{m+1}} \, G^{m+1} \biggr)
   \,-\,
   G^{m+1} \biggl( {H'(\zeta,w) \over (m+1) \zeta^{m+1}}
                   \,-\,
                   {H(\zeta,w) \over \zeta^{m+2}}  \biggr)
   \nonumber \\
\end{eqnarray}
where the prime denotes $\partial/\partial\zeta$.
Since the total derivative gives zero when integrated around a closed contour,
we have
\begin{equation}
   H(\varphi(w),w)  \;=\;
   \sum_{m=0}^\infty
   {1 \over 2\pi i} \oint\limits_{|\zeta| = \rho}
   \left[ {H(\zeta,w) \over \zeta^{m+1}} G^m
          \,+\,
          {H'(\zeta,w) \over (m+1) \zeta^{m+1}} G^{m+1}
          \,-\,
          {H(\zeta,w) \over \zeta^{m+2}} G^{m+1}
   \right] \, d\zeta
   \;.
 \label{eq.after.intbyparts}
\end{equation}
Now the first and third terms in brackets cancel when summed over $m$,
except for the first term at $m=0$, which gives simply
$(1/2 \pi i) \oint [H(\zeta,w)/\zeta] \, d\zeta = H(0,w)$.
Hence
\begin{eqnarray}
   H(\varphi(w),w)
   & = &
   H(0,w) \,+\, \sum_{m=0}^\infty
   {1 \over 2\pi i} \oint
          {H'(\zeta,w) \over (m+1) \zeta^{m+1}} G(\zeta,w)^{m+1}  \, d\zeta
   \nonumber \\[2mm]
   & = &
   H(0,w) \,+\, \sum_{m=1}^\infty
   {1 \over 2\pi i} \oint
          {H'(\zeta,w) \over m \zeta^{m}} G(\zeta,w)^{m}  \, d\zeta
   \nonumber \\[2mm]
   & = &
   H(0,w) \,+\, \sum_{m=1}^\infty {1 \over m} \,
        [\zeta^{m-1}]  H'(\zeta,w) G(\zeta,w)^{m}
   \;.
 \label{eq.formula.hvarphi}
\end{eqnarray}
This proves the fundamental formula \reff{eq.star4}/\reff{eq.star4bis},
in which the sum is absolutely convergent on $V$,
locally uniformly on $V$.

It follows from this formula [taking $H(z,w)=z$]
that if $V \subset \C^M$ and $G$ is analytic,
then $\varphi$ is analytic as well;
and if also $H$ is analytic, then so is $w \mapsto H(\varphi(w),w)$.

Finally, if $V$ is a polydisc
(or more generally a complete Reinhardt domain)
centered at $0 \in \C^M$,
and $G$ and $H$ are analytic,
then the analytic function $H(\varphi(w),w)$ is given in $V$
by a convergent Taylor series.
The coefficients of this Taylor series are given by
\reff{eq.calpha}/\reff{eq.calpha.a}
because the absolutely convergent function series
\reff{eq.star4bis}/\reff{eq.star4abis}
can be differentiated term-by-term.
\qed

\begin{remark}
\rm
The following alternative calculation
provides a slightly slicker proof of Theorem~\ref{thm1}:
start from \reff{eq.cauchy} and write
\begin{eqnarray}
   H(\zeta,w) \, {1 - {\partial G \over \partial\zeta}(\zeta,w)
                  \over
                  \zeta - G(\zeta,w)
                 }
   & = &
   H(\zeta,w) \, {\partial \over \partial\zeta} \,
                 \log\bigl[ \zeta - G(\zeta,w) \bigr]
       \nonumber \\[2mm]
   & = &
   H(\zeta,w) \left(  {1 \over \zeta} \,+\,
                      {\partial \over \partial\zeta} \,
                      \log\biggl[ 1 - {G(\zeta,w) \over \zeta} \biggr]
              \right)
   \,,
 \label{eq.remark.altproof}
\end{eqnarray}
where hypothesis (c) and Taylor expansion guarantee that
the function $\log[1 - G(\zeta,w)/\zeta]$
is well-defined and single-valued on the circle $|\zeta| = \rho$;
furthermore, a simple compactness argument extends this
to some annulus $\rho_1 < |\zeta| \le \rho$
(locally uniformly in $w$).
Now, the first term in \reff{eq.remark.altproof}, when integrated,
yields $H(0,w)$.
To handle the second term, let us integrate by parts, Taylor-expand the log,
and then extract the residue from the resulting Laurent series:  we get
\begin{eqnarray}
   - \, {1 \over 2\pi i} \oint
          H'(\zeta,w) \log\biggl[ 1 - {G(\zeta,w) \over \zeta} \biggr] \, d\zeta
   & = &
   \sum_{m=1}^\infty {1 \over m} [\zeta^{-1}] H'(\zeta,w)
                          \biggl( {G(\zeta,w) \over \zeta} \biggr)^m
   \nonumber \\[2mm]
   & = &
   \sum_{m=1}^\infty {1 \over m} [\zeta^{m-1}] H'(\zeta,w) G(\zeta,w)^m
   \;,
   \qquad\qquad
 \label{eq.log.altproof}
\end{eqnarray}
which is \reff{eq.star4}/\reff{eq.star4bis}.
A similar argument without integration by parts yields
the alternate formula \reff{eq.star4a}/\reff{eq.star4abis}.
I~thank Alex Eremenko for helpful comments concerning this proof.
\qed
\end{remark}

\begin{remark}
  \label{remark.gessel}
\rm
Formula \reff{eq.star4}/\reff{eq.star4bis} can alternatively be deduced
from the standard Lagrange inversion formula \reff{eq.star2}
by an argument due to Ira Gessel \cite[p.~148]{Stanley_99}:
Introduce a new parameter $t \in \C$, and study the equation $z = t G(z,w)$
with solution $z = \Phi(w,t)$.
Applying \reff{eq.star2} with $w$ fixed and $t$ as the variable, we obtain
\begin{equation}
   H(\Phi(w,t),w)  \;=\;
   H(0,w) \,+\,
   \sum_{m=1}^\infty {t^m \over m} \,
      [\zeta^{m-1}] {\partial H(\zeta,w) \over \partial\zeta} \, G(\zeta,w)^m
   \;.
\end{equation}
Setting $t=1$ yields \reff{eq.star4}/\reff{eq.star4bis}.
[Hypothesis (c) guarantees that
 $\sup\limits_{|\zeta| = \rho} |G(\zeta,w)/\zeta| \equiv C < 1$,
 so that the Lagrange series is convergent for $|t| < 1/C$.]
An analogous argument starting from \reff{eq.star2a} yields
the alternate formula \reff{eq.star4a}/\reff{eq.star4abis}.
\qed
\end{remark}

\begin{remark}
\rm
For some purposes the disc $\D_\rho$ can be replaced an arbitrary
domain (connected open set) $D \subset \C$.
Hypothesis (c) is then replaced by the assumption that the image
$G(D,w)$ is relatively compact in $D$ for all $w \in V$.
Under this hypothesis, the map $z \mapsto G(z,w)$
is a strict contraction in the Poincar\'e metric \cite{Krantz_08} on $G(D,w)$,
i.e.\ satisfies
$d_{\rm Poin}(G(z_1,w),\, G(z_2,w)) \le \kappa d_{\rm Poin}(z_1,z_2)$
for some $\kappa < 1$ (locally uniformly in $w$).
It then follows from the contraction-mapping fixed-point theorem
that there is a unique fixed point $\varphi(w)$,
which moreover can be obtained by iteration starting at any point of $D$.
That is, if we define
\begin{subeqnarray}
   G_0(\zeta,w)  & = &  \zeta   \\[1mm]
   G_{n+1}(\zeta,w)  & = &  G(G_n(\zeta,w),w)
\end{subeqnarray}
then $G_n(\zeta,w) \to \varphi(w)$,
uniformly for $\zeta \in D$ and locally uniformly for $w \in V$.
It would be interesting to know whether
this can be used to provide a function series
analogous to \reff{eq.star3}--\reff{eq.star4a}
based on the Taylor coefficients of $G$ at
an arbitrary point $\zeta_0 \in D$.
An analogous argument works for domains $D \subset \C^N$,
using the Kobayashi metric
\cite{Lang_87,Jarnicki_93,Kobayashi_98,Kobayashi_05,Krantz_08},
provided that some iterate $G^n(D,w)$ is Kobayashi-hyperbolic.
I~thank Alex Eremenko for suggesting the use of the Poincar\'e
and Kobayashi metrics.
\qed
\end{remark}

\begin{question}
  \label{question.pole}
\rm
Can Theorem~\ref{thm1} be generalized to allow $H(z,w)$ to have a pole
at $z=0$?  Please note that by linearity it suffices to consider
$H(z,w) = h(w) z^{-k}$ ($k \ge 1$);
and since $h(w)$ just acts as an overall prefactor,
it suffices to consider simply $H(z,w) = z^{-k}$.
Of course we will somehow have to restrict attention
to the subset of $V$ where $\varphi(w) \neq 0$;
and some hypothesis will be needed to guarantee that this subset is nonempty,
i.e.\ that $\varphi$ is not identically zero.
\qed
\end{question}

\section{Implicit function formula: Formal-power-series version}

In this section we shall consider $G(z,w)$ to be a formal power series
in indeterminates $z$ and $w = (w_i)_{i \in I}$,
where $I$ is an arbitrary finite or infinite index set.
The coefficients in this formal power series may belong
to an arbitrary commutative ring-with-identity-element $R$.
(For some purposes we will want to assume further that
the coefficient ring $R$ contains the rationals as a subring.
 In applications, $R$ will usually be a field of characteristic 0
 --- e.g.\ the rationals $\Q$, the reals $\R$, or the complex numbers $\C$ ---
 or a ring of polynomials or formal power series over such a field.)
See \cite{Cartan_63,Niven_69,Henrici_74,Hungerford_74,Goulden_83}
for basic facts about formal power series.
We recall that $R[[w]]$ denotes the ring of formal power series
in the indeterminates $w = (w_i)_{i \in I}$ with coefficients in $R$.

We begin with a well-known proposition asserting the existence and uniqueness
of the formal power series $\varphi(w)$ solving the equation $z = G(z,w)$,
or equivalently $F(z,w)=0$.
Here $R$ can be an arbitrary commutative ring with identity element;
it need not contain the rationals or even be of characteristic 0.
(For instance, the ring $\Z_n$ of integers modulo $n$ is allowed.)


\begin{proposition}[Implicit function theorem for formal power series]
   \label{prop.existence.uniqueness}
\mbox{}
\begin{itemize}
   \item[(a)] Let $R$ be a commutative ring with identity element.
Let $F(z,w)$ be a formal power series in indeterminates $z$ and
$w = (w_i)_{i \in I}$, with coefficients in $R$;
suppose further that $F(0,0) = 0$
and that $(\partial F/\partial z)(0,0)$ is invertible in the ring $R$.
Then there exists a unique formal power series $\varphi(w)$
with zero constant term satisfying $F(\varphi(w),w) = 0$.
   \item[(b)] Let $R$ be a commutative ring with identity element.
Let $G(z,w)$ be a formal power series in indeterminates $z$ and
$w = (w_i)_{i \in I}$, with coefficients in $R$;
suppose further that $G(0,0) = 0$
and that $1 - (\partial G/\partial z)(0,0)$ is invertible in the ring $R$.
Then there exists a unique formal power series $\varphi(w)$
with zero constant term satisfying $\varphi(w) = G(\varphi(w),w)$.
\end{itemize}
\end{proposition}

If the ring $R$ is a field,
then the hypothesis that $(\partial F/\partial z)(0,0)$
be invertible in $R$ means simply that $(\partial F/\partial z)(0,0) \neq 0$.
If $R$ is a ring of formal power series over a field,
then this hypothesis
means that the constant term of $(\partial F/\partial z)(0,0)$ is $\neq 0$.
Analogous statements apply to part (b),
with $(\partial G/\partial z)(0,0) \neq 1$.

\proofof{Proposition~\ref{prop.existence.uniqueness}}
It suffices to prove either (a) or (b), since they are equivalent
under the substitution $F(z,w) = z-G(z,w)$.
We shall prove (b).  Let us write
\begin{equation}
   G(z,w)  \;=\;  \sum_{k=0}^\infty g_k(w) \, z^k
\end{equation}
with each $g_k \in R[[w]]$;
by hypothesis $g_0(0) = 0$, and $1 - g_1(0)$ is invertible in the ring $R$.
The equation $\varphi(w) = G(\varphi(w),w)$ can now be written as
\begin{equation}
   \varphi(w)  \;=\; \sum_{k=0}^\infty g_k(w) \, \varphi(w)^k
\end{equation}
or equivalently
\begin{equation}
   \varphi(w)
   \;=\;
   [1 - g_1(0)]^{-1} 
   \left[ g_0(w) \:+\: [g_1(w)-g_1(0)] \, \varphi(w) \:+\:
          \sum_{k=2}^\infty g_k(w) \, \varphi(w)^k
   \right]   \,.
 \label{eq.prop.existence.uniqueness.proof}
\end{equation}
Since the series $g_0(w)$, $g_1(w)-g_1(0)$ and $\varphi(w)$
all have zero constant term,
we see from \reff{eq.prop.existence.uniqueness.proof}
that the coefficients in
$\varphi(w) = \sum_{|\balpha| \ge 1} b_\balpha w^\balpha$
can be uniquely determined by induction on $|\balpha|$.
Conversely, the unique solution of this system of equations
necessarily solves $\varphi(w) = G(\varphi(w),w)$.
\qed

A multidimensional generalization of
Proposition~\ref{prop.existence.uniqueness},
in which $z$ is replaced by a vector of indeterminates $(z_1,\ldots,z_N)$,
can be found in Bourbaki \cite[p.~A.IV.37]{Bourbaki_90}.
As one might expect, the hypothesis is that the Jacobian determinant
$\det(\partial F/\partial z)(0,0)$ is invertible in $R$.

We can carry this argument further and provide an explicit formula
for $\varphi(w)$.
Let us begin with what appears to be a special case,
but in fact contains the general result:
namely, let us take $G(z,\g) = \sum_{n=0}^\infty g_n z^n$
where $\g = (g_n)_{n=0}^\infty$ are indeterminates.
We then have the following ``universal'' version of the
Lagrange inversion formula \cite[Theorem~6.2]{Gessel-Stanley}:

\begin{proposition}[Universal Lagrange inversion formula]
   \label{prop.universal_Lagrange}
Let $\g = (g_n)_{n=0}^\infty$ be indeterminates.
There is a unique formal power series $\varphi \in \Z[[\g]]$
with zero constant term satisfying
$\varphi(\g) = \sum\limits_{n=0}^\infty g_n \varphi(\g)^n$,
and its coefficients are given explicitly by
\begin{equation}
   [g_0^{k_0} g_1^{k_1} g_2^{k_2} \cdots] \, \varphi(\g)^\ell
   \;=\;
   \begin{cases}
      \ell \: {\textstyle \Bigl( \, \sum\limits_{n=0}^\infty k_n \,-\, 1 \Bigr)!
               \over
               \textstyle \prod\limits_{n=0}^\infty k_n!
              }
         & \textrm{if $\sum\limits_{n=0}^\infty (n-1) k_n = -\ell$}  \\[2mm]
      0  & \textrm{otherwise}
    \end{cases}
 \label{eq.prop.universal_Lagrange}
\end{equation}
for all integers $\ell \ge 1$.
In particular, for each pair $(k_0,k_1)$
there are {\em finitely many}\/ $(k_2,k_3,\ldots)$
for which \reff{eq.prop.universal_Lagrange} is nonzero,
so that $\varphi(\g)^\ell$ is a formal power series in $g_0,g_1$
whose coefficients are {\em polynomials}\/ in $g_2,g_3,\ldots$:
that is, $\varphi \in \Z[g_2,g_3,\ldots][[g_0,g_1]]$.

We can also write the formula
\begin{equation}
   \varphi(\g)^\ell
   \;=\;
   \sum_{m=1}^\infty  [\zeta^{m-\ell}]
       \left( \sum\limits_{n=0}^\infty g_n \zeta^n \right)^{\! m-1}
       \!
       \left( \sum\limits_{n=0}^\infty (1-n) g_n \zeta^n \right)
   \;,
 \label{eq.prop.universal_Lagrange.2alt}
\end{equation}
and in the ring $\Q[[\g]]$ we can write
\begin{equation}
   \varphi(\g)^\ell
   \;=\;
   \sum_{m=1}^\infty {\ell \over m} \, [\zeta^{m-\ell}]
       \left( \sum\limits_{n=0}^\infty g_n \zeta^n \right)^{\! m}
   \;.
 \label{eq.prop.universal_Lagrange.2}
\end{equation}
\end{proposition}

\begin{figure}[t]
\setlength{\unitlength}{1.5cm}
\begin{center}
\begin{picture}(9,4)
\put(1,3){\circle*{0.15}}
\put(1,2.65){\circle{0.7}}
\put(0.78,2.59){$\footnotesize \varphi(\g)$}
\put(1.9,2.65){$=$}
\put(3,3){\circle*{0.15}}
\put(2.9,3.3){$g_0$}
\put(3.62,2.65){$+$}
\put(4.5,3){\circle*{0.15}}
\put(4.4,3.3){$g_1$}
\put(4.5,3){\line(0,-1){1}}
\put(4.5,2){\circle*{0.15}}
\put(4.5,1.65){\circle{0.7}}
\put(4.28,1.59){$\footnotesize \varphi(\g)$}
\put(5.2,2.65){$+$}
\put(6.4,3){\circle*{0.15}}
\put(6.3,3.3){$g_2$}
\put(6.4,3){\line(-1,-2){0.5}}
\put(6.4,3){\line(1,-2){0.5}}
\put(5.9,2){\circle*{0.15}}
\put(6.9,2){\circle*{0.15}}
\put(5.9,1.65){\circle{0.7}}
\put(6.9,1.65){\circle{0.7}}
\put(5.68,1.59){$\footnotesize \varphi(\g)$}
\put(6.68,1.59){$\footnotesize \varphi(\g)$}
\put(7.4,2.65){$+\;\;\ldots$}
\end{picture}
\end{center}
\vspace*{-2cm}
\caption{
   The functional equation
   $\varphi(\g) = \sum_{n=0}^\infty g_n \varphi(\g)^n$
   for the generating function of unlabeled plane trees.
}
 \label{fig1}
\end{figure}
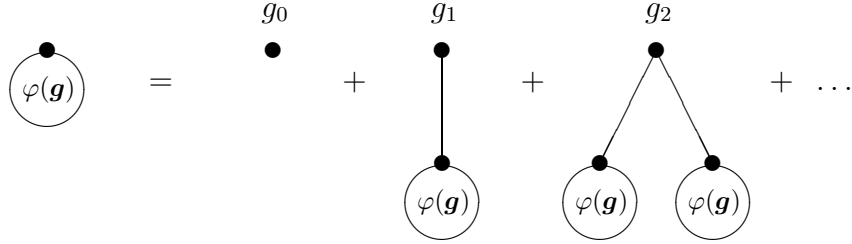

\proof
The functional equation
$\varphi(\g) = \sum_{n=0}^\infty g_n \varphi(\g)^n$
is the same equation as is satisfied by
the (ordinary) generating function for {\em unlabeled plane trees}\/
(i.e., rooted trees in which the vertices are unlabeled
 but the subtrees at each vertex are linearly ordered),
in which a vertex having $n$ children gets a weight $g_n$,
and the weight of a tree is the product of its vertex weights
(see Figure~\ref{fig1}).\footnote{
   The key fact here is that if $\Sigma$ is any sum of terms,
 then a term in the expansion of $\Sigma^n = \Sigma \Sigma \cdots \Sigma$
 is obtained by choosing, {\em in order}\/,
 a term of $\Sigma$ for the first factor,
 a term of $\Sigma$ for the second factor, etc.\
 This is why one obtains trees in which the subtrees at each vertex
 are linearly ordered.
}
Since the solution of this equation is unique
by Proposition~\ref{prop.existence.uniqueness}(b),
it follows that $\varphi(\g)$ is
the generating function for unlabeled plane trees with this weighting.
More generally, $\varphi(\g)^\ell$
is the generating function for
{\em unlabeled plane forests}\/ with $\ell$ components,
with the same weighting.\footnote{
   An {\em unlabeled plane forest}\/ is a forest of rooted trees
   with unlabeled vertices
   in which the subtrees at each vertex are linearly ordered
   and the components of the forest (or equivalently their roots)
   are also linearly ordered.
   The reasoning in the preceding footnote explains why
   $\varphi(\g)^\ell$ gives rise to forests in which
   the components are linearly ordered.
}
It is a well-known (though nontrivial) combinatorial fact
\cite[Theorem~5.3.10]{Stanley_99}
that the number of unlabeled plane forests with $\ell$ components
having type sequence $(k_0,k_1,\ldots)$
[i.e., in which there are $k_n$ vertices having $n$ children,
 for each $n \ge 0$]
is given precisely by \reff{eq.prop.universal_Lagrange}.

Since the constraint in \reff{eq.prop.universal_Lagrange}
can be written as $k_0 = \ell + \sum_{n=2}^\infty (n-1) k_n$,
one sees immediately that for each pair $(k_0,k_1)$
there are finitely many $(k_2,k_3,\ldots)$
for which \reff{eq.prop.universal_Lagrange} is nonzero.

To prove \reff{eq.prop.universal_Lagrange.2alt},
let us expand out the summand on the right-hand side,
choosing $n=N$ in the last factor:
we get
\begin{eqnarray}
   & & \!\!\!\!\!
   [\zeta^{m-\ell}]
       \left( \sum\limits_{n=0}^\infty g_n \zeta^n \right)^{\! m-1}
       \!
       \left( \sum\limits_{n=0}^\infty (1-n) g_n \zeta^n \right)
   \nonumber \\
   & & \;=\;
   \sum\limits_{N=0}^\infty
   [\zeta^{m-\ell}]
   \!
         \sum_{\substack{
                  k_0,k_1,k_2,\ldots \ge 0 \\[1mm]
                  \sum k_n = m
              }}
         {m-1 \choose k_0,\ldots,k_{N-1}, k_N-1, k_{N+1},\ldots} \,
         (1-N)
         \left( \prod_{n=0}^\infty g_n^{k_n} \right)
         \zeta^{\sum n k_n}
   \nonumber \\
\end{eqnarray}
where ${\displaystyle {m-1 \choose k_0,\ldots,k_{N-1}, k_N-1, k_{N+1},\ldots}}
  = (m-1)! \, k_N \Big/ \prod\limits_{n=0}^\infty k_n!$
is a multinomial coefficient.
The constraints tell us that $\sum_{n=0}^\infty k_n = m$
and $\sum_{n=0}^\infty n k_n = m-\ell$,
so that $\sum_{N=0}^\infty (1-N) k_N = \ell$.
This is precisely \reff{eq.prop.universal_Lagrange}.

If work over the rationals rather than the integers,
things become slightly easier:
expanding out the summand in \reff{eq.prop.universal_Lagrange.2}, we obtain
\begin{equation}
   {\ell \over m} \, [\zeta^{m-\ell}]
       \left( \sum\limits_{n=0}^\infty g_n \zeta^n \right)^{\! m}
   \;=\;
   {\ell \over m} \, [\zeta^{m-\ell}]
         \sum_{\substack{
                  k_0,k_1,k_2,\ldots \ge 0 \\[1mm]
                  \sum k_n = m
              }}
         {m \choose k_0,k_1,k_2,\ldots} \,
         \left( \prod_{n=0}^\infty g_n^{k_n} \right)
         \zeta^{\sum n k_n}
   \;,
\end{equation}
which again matches \reff{eq.prop.universal_Lagrange}.
\qed

\begin{remark}
\rm
To see easily (i.e., without the full combinatorial interpretation)
that the numbers \reff{eq.prop.universal_Lagrange} are indeed integers,
it suffices to note that
\begin{equation}
   k_i \:  {\textstyle \bigl( \, \sum k_n \,-\, 1 \bigr)!
            \over
            \textstyle \prod k_n!
           }
   \;\,=\;\,
   { \textstyle \sum k_n \,-\, 1
     \choose
     \textstyle k_1,\ldots,k_{i-1}, k_i-1, k_{i+1},\ldots
   }
\end{equation}
is a multinomial coefficient and hence an integer,
and that $\ell = \sum_{i=0}^\infty (1-i) k_i$
by virtue of the constraint.
I~thank Richard Stanley and Ira Gessel (independently) for this observation.
\qed
\end{remark}

\begin{remark}
\rm
Let us stress that the proof of the counting result
given in \cite[Theorem~5.3.10]{Stanley_99} is purely combinatorial;
it is based on a bijection between unlabeled plane forests
and a certain class of words on a finite alphabet.
We shall use it to deduce the implicit function formula
(and in particular the Lagrange inversion formula)
as a simple corollary, following the second proof of
\cite[Theorem~5.4.2]{Stanley_99}.
This approach to the Lagrange inversion formula 
goes back to Raney \cite{Raney_60}
and was later simplified by other authors \cite{Schutzenberger_71,Lothaire_83}.
On the other hand, a much easier (though perhaps less enlightening)
way of obtaining this counting result is to first prove
the Lagrange inversion formula
(e.g.\ by the algebraic argument given in the first proof of
 \cite[Theorem~5.4.2]{Stanley_99})
and then use it to obtain the enumeration of plane trees or forests
as a straightforward application
\cite{Tutte_64} \cite[section~2.7.7]{Goulden_83}.
\qed
\end{remark}

We can now deduce the general implicit function formula
as an easy corollary of Proposition~\ref{prop.universal_Lagrange}.
The key point is that $\varphi(\g)^\ell$ is a formal power series in $g_0,g_1$
whose coefficients are {\em polynomials}\/ in $g_2,g_3,\ldots$:
therefore, in the identity
$\varphi(\g) = \sum_{n=0}^\infty g_n \varphi(\g)^n$
we can make the substitutions $g_n \leftarrow g_n(w)$,
where the $g_n(w)$ are formal power series in an arbitrary collection
$w = (w_i)_{i \in I}$ of indeterminates,
provided that $g_0(w)$ and $g_1(w)$ have zero constant term;
the series $g_n(w)$ for $n \ge 2$ are unrestricted.
This yields a formal-power-series version of \reff{eq.star4a}:

\begin{theorem}[Implicit function formula --- formal-power-series version \#1]
 \label{thm2.first}
Let $R$ be a commutative ring with identity element.
Let $G(z,w) = \sum_{n=0}^\infty g_n(w) \, z^n$
be a formal power series in indeterminates $z$ and $w = (w_i)_{i \in I}$,
with coefficients in $R$,
satisfying $G(0,0) = 0$ and $(\partial G/\partial z)(0,0) = 0$
[i.e., $g_0(0) = 0$ and $g_1(0) = 0$].
Then the unique formal power series $\varphi(w)$
with zero constant term satisfying $\varphi(w) = G(\varphi(w),w)$
is given explicitly by
\begin{equation}
   \varphi(w)^\ell  \;=\;
   \sum_{k_0,k_1,k_2,\ldots \ge 0}
   c_\ell(k_0,k_1,k_2,\ldots) \prod_{n=0}^\infty g_n(w)^{k_n}
 \label{eq.thm2.first}
\end{equation}
for all integers $\ell \ge 1$,
where $c_\ell(k_0,k_1,k_2,\ldots)$ is given by
\reff{eq.prop.universal_Lagrange}.
Its coefficients are given by the finite sums
\begin{equation}
   [w^\balpha] \varphi(w)  \;=\;
   \sum_{m=1}^{2|\balpha|}
      [\zeta^{m} w^\balpha] \Bigl[ G(\zeta,w)^m \:-\:
                     \zeta \, {\partial G(\zeta,w) \over \partial\zeta} \,
                     G(\zeta,w)^{m-1}
                              \Bigr]
   \;.
 \label{eq.thm2.first.eq1}
\end{equation}
More generally, if $H(z,w)$ is any formal power series, we have
\begin{eqnarray}
   & & \!\!\!\!\!
   [w^\balpha] H(\varphi(w),w)
        \nonumber \\
   & &
   \;=\;
   [w^\balpha] H(0,w)  \,+\,
   \sum_{m=1}^{2|\balpha|}
      [\zeta^{m} w^\balpha]
        H(\zeta,w) \, \Bigl[ G(\zeta,w)^m \:-\:
                     \zeta \, {\partial G(\zeta,w) \over \partial\zeta} \,
                     G(\zeta,w)^{m-1}
                              \Bigr]
   \;.
   \nonumber \\ \label{eq.thm2.first.eq2}
\end{eqnarray}
\end{theorem}

\proof
As just observed,
substituting $g_n \leftarrow g_n(w)$ in \reff{eq.prop.universal_Lagrange}
proves \reff{eq.thm2.first}.

To prove \reff{eq.thm2.first.eq2},
suppose first that $H(z,w) = z^\ell$.
If $\ell=0$, then obviously $\varphi(w)^\ell = 1$.
If $\ell \ge 1$, we can substitute $g_n \leftarrow g_n(w)$
in \reff{eq.prop.universal_Lagrange.2alt} to obtain
\begin{equation}
   \varphi(w)^\ell
   \;=\;
   \sum_{m=1}^\infty [\zeta^m] \, \zeta^{\ell}
           \, \Bigl[ G(\zeta,w)^m \:-\:
                     \zeta \, {\partial G(\zeta,w) \over \partial\zeta} \,
                     G(\zeta,w)^{m-1}
              \Bigr]
   \;.
 \label{eq.proof.thm2.star1.alt}
\end{equation}
The general case of $H(z,w) = \sum_{\ell=0}^\infty h_\ell(w) z^\ell$
is obtained from these special cases by linearity:  we have
\begin{equation}
   H(\varphi(w),w)
   \;=\;
   H(0,w) \,+\,
   \sum_{m=1}^\infty
      [\zeta^m] H(\zeta,w)  \Bigl[ G(\zeta,w)^m \:-\:
                     \zeta \, {\partial G(\zeta,w) \over \partial\zeta} \,
                     G(\zeta,w)^{m-1}
                              \Bigr]
   \;.
\end{equation}
Now take the coefficient of $w^\balpha$ on both sides,
and observe that contributions come only from $m \le 2|\balpha|$
(see Remark~\ref{remark.calpha}).
\qed

The implicit function formulae given in Theorem~\ref{thm2.first}
--- in particular, the variant Lagrange formula
\reff{eq.thm2.first.eq1}/\reff{eq.thm2.first.eq2} ---
are valid in an arbitrary commutative ring-with-identity-element $R$,
even if $R$ is not of characteristic 0,
because the numerical coefficients arising in
\reff{eq.thm2.first}--\reff{eq.thm2.first.eq2} are all integers.
On the other hand, if we are willing to assume that the ring $R$
contains the rationals as a subring,
then we can deduce the slightly more convenient
explicit formulae \reff{eq.star3}/\reff{eq.star4} for $\varphi(w)$.
An argument completely analogous to that leading to \reff{eq.thm2.first.eq2},
but based on \reff{eq.prop.universal_Lagrange.2}
instead of \reff{eq.prop.universal_Lagrange.2alt}, proves:

\begin{theorem}[Implicit function formula --- formal-power-series version \#2]
 \label{thm2}
Let $R$ be a commutative ring containing the rationals as a subring.
Let $G(z,w)$ be a formal power series
in indeterminates $z$ and $w = (w_i)_{i \in I}$,
with coefficients in $R$,
satisfying $G(0,0) = 0$ and $(\partial G/\partial z)(0,0) = 0$.
Then there exists a unique formal power series $\varphi(w)$
with zero constant term satisfying $\varphi(w) = G(\varphi(w),w)$,
and its coefficients are given by the finite sums
\begin{equation}
   [w^\balpha] \varphi(w)  \;=\;
   \sum_{m=1}^{2 |\balpha| - 1}
      {1 \over m} \, [\zeta^{m-1} w^\balpha]  G(\zeta,w)^m
   \;.
 \label{eq.thm2.eq1}
\end{equation}
More generally, if $H(z,w)$ is any formal power series, we have
\begin{equation}
   [w^\balpha] H(\varphi(w),w)  \;=\;
   [w^\balpha] H(0,w)  \,+\,
   \sum_{m=1}^{2 |\balpha| - 1}
      {1 \over m} \, [\zeta^{m-1} w^\balpha]
        {\partial H(\zeta,w) \over \partial\zeta} \, G(\zeta,w)^m
   \;.
 \label{eq.thm2.eq2}
\end{equation}
\end{theorem}



It is instructive to give two alternate proofs of Theorem~\ref{thm2}:
one deducing the result as a corollary of the usual Lagrange inversion
formula for formal power series, and one deducing it from 
our analytic version of the implicit function formula (Theorem~\ref{thm1}).

\secondproofof{Theorem~\ref{thm2}}
Proposition~\ref{prop.existence.uniqueness}(b)
gives the existence and uniqueness of $\varphi(w)$.
[The proof of \reff{eq.thm2.eq1}/\reff{eq.thm2.eq2} to be given next
 will provide an alternate proof of uniqueness.]
Now we deduce \reff{eq.thm2.eq2} by using the
formal-power-series version of Gessel's argument \cite[p.~148]{Stanley_99}
mentioned in Remark~\ref{remark.gessel}:
Introduce a new indeterminate $t$, and study the equation $z = t G(z,w)$
with solution $z = \Phi(w,t)$.
Using the Lagrange inversion formula \reff{eq.star2}
for formal power series in the indeterminate $t$
and with coefficients in the ring $R[[w]]$, we obtain\footnote{
   The Lagrange inversion formula for formal power series
   is most commonly stated for series with coefficients in a
   field of characteristic 0 \cite[Theorem~5.4.2]{Stanley_99},
   but in the form \reff{eq.star1}/\reff{eq.star2}
   it also holds for series with coefficients in an
   arbitrary commutative ring containing the rationals.
   See e.g.\ \cite[Theorem~1.2.4]{Goulden_83}.
}
\begin{equation}
   H(\Phi(w,t),w)  \;=\;
   H(0,w) \,+\,
   \sum_{m=1}^\infty {t^m \over m} \,
      [\zeta^{m-1}] {\partial H(\zeta,w) \over \partial\zeta} \, G(\zeta,w)^m
\end{equation}
and hence
\begin{equation}
   [w^\balpha] H(\Phi(w,t),w)  \;=\;
   [w^\balpha] H(0,w) \,+\,
   \sum_{m=1}^\infty {t^m \over m} \,
      [\zeta^{m-1} w^\balpha] {\partial H(\zeta,w) \over \partial\zeta} \,
          G(\zeta,w)^m
 \label{eq.thm2.proof}
\end{equation}
as equalities between formal power series in $t$.
But by the hypothesis on $G$, the only nonzero contributions
to the sum on the right-hand side of \reff{eq.thm2.proof}
come from $m \le 2 |\balpha| - 1$ (see Remark~\ref{remark.calpha}),
so each side of \reff{eq.thm2.proof} is in fact a {\em polynomial}\/
in $t$ (of degree at most $2 |\balpha| - 1$).
So we can evaluate it at any chosen $t \in R$,
in particular at $t=1$.
This proves \reff{eq.thm2.eq2}.
\qed

\thirdproofof{Theorem~\ref{thm2}}
We see from the proof of Proposition~\ref{prop.existence.uniqueness}
that $\varphi(w)^\ell$ will be given by a universal formula
of the form \reff{eq.thm2.first} with nonnegative integer coefficients
$c_\ell(k_0,k_1,k_2,\ldots)$;  it remains only to find these coefficients.
To do this, it suffices to consider the case $R = \C$
with a single indeterminate $w$.\footnote{
   Here it is important that $\C$ is of characteristic 0,
   in order to avoid losing information about the integers
   $c_\ell(k_0,k_1,k_2,\ldots)$.
}
It furthermore suffices to consider the cases in which $G(z,w)$
is a {\em polynomial}\/ in $z$ and $w$ (of arbitrarily high degree),
since $[w^n] \varphi(w)$ depends only on the $[z^j w^k] G(z,w)$
with $j,k \le n$.  But in this case we can apply
the analytic version of the implicit function formula (Theorem~\ref{thm1}).
\qed


%
%

\begin{question}
  \label{question.pole.powerseries}
\rm
Can Theorems~\ref{thm2.first} and \ref{thm2} be generalized to allow $H(z,w)$
to be a {\em Laurent}\/ series in $z$,
at least when $w$ is a single indeterminate
and $b_{01} \equiv (\partial G/\partial w)(0,0)$ is
invertible in the ring $R$?
(Or slightly more restrictively,
 when $R$ is a field of characteristic zero and $b_{01} \neq 0$?)
Please note that by linearity it suffices to consider
$H(z,w) = h(w) z^{-k}$ ($k \ge 1$);
and since $h(w)$ just acts as an overall prefactor,
it suffices to consider simply $H(z,w) = z^{-k}$.
One might try imitating \cite[first proof of Theorem~5.4.2]{Stanley_99},
possibly combined with the Gessel idea $z = t G(z,w)$.
\qed
\end{question}


In Theorems~\ref{thm2.first} and \ref{thm2},
we have for simplicity assumed that
$b_{10} \equiv (\partial G/\partial z)(0,0) = g_1(0)$ is {\em zero}\/.
This seems more restrictive than Proposition~\ref{prop.existence.uniqueness},
where it was assumed only that $1-b_{10}$ is invertible in the ring $R$,
but there is in fact no real loss of generality here.
For if $1-b_{10}$ is invertible, then the equation $z=G(z,w)$
is equivalent to $z= \widetilde{G}(z,w)$, where
\begin{equation}
   \widetilde{G}(z,w)  \;=\;
   (1-b_{10})^{-1} \, [G(z,w) - b_{10}z]
\end{equation}
satisfies $(\partial \widetilde{G}/\partial z)(0,0) = 0$.
We can therefore apply Theorems~\ref{thm2.first} and \ref{thm2} with
\begin{equation}
   \widetilde{g}_n(w)  \;=\;
   [1 - g_1(0)]^{-1} \,\times\,
   \left\{ \! \begin{array}{ll}
                 g_n(w)          &   \text{if $n \neq 1$}   \\[1mm]
                 g_1(w)-g_1(0)   &   \text{if $n = 1$}
              \end{array}
   \! \right\}
   \;.
\end{equation}

On the other hand, if $R = \R$ or $\C$ and $|b_{10}| < 1$,
we can avoid this preliminary transformation if we prefer:
the power-series coefficients will then be given by
absolutely convergent infinite sums, which are nothing other than
the finite sums based on $\widetilde{G}$ in which each factor
$(1-b_{10})^{-m}$ has been expanded out as
$\sum\limits_{k=0}^\infty {\displaystyle {m+k-1 \choose k}} \, b_{10}^k$.
We therefore have:

\begin{theorem}[Implicit function formula --- formal-power-series version \#3]
 \label{thm2bis}
Let $R$ be either the field $\R$ of real numbers
or the field $\C$ of complex numbers.
Let $G(z,w)$ be a formal power series
in indeterminates $z$ and $w = (w_i)_{i \in I}$,
with coefficients in $R$, satisfying $G(0,0) = 0$
and $|(\partial G/\partial z)(0,0)| < 1$.
Then there exists a unique formal power series $\varphi(w)$
with zero constant term satisfying $\varphi(w) = G(\varphi(w),w)$,
and its coefficients are given by the absolutely convergent sums
\begin{equation}
   [w^\balpha] \varphi(w)  \;=\;
   \sum_{m=1}^{\infty}
      {1 \over m} \, [\zeta^{m-1} w^\balpha]  G(\zeta,w)^m
   \;.
 \label{eq.thm2bis.eq1}
\end{equation}
More generally, if $H(z,w)$ is any formal power series, we have
\begin{equation}
   [w^\balpha] H(\varphi(w),w)  \;=\;
   [w^\balpha] H(0,w)  \,+\,
   \sum_{m=1}^{\infty}
      {1 \over m} \, [\zeta^{m-1} w^\balpha]
        {\partial H(\zeta,w) \over \partial\zeta} \, G(\zeta,w)^m
   \;.
 \label{eq.thm2bis.eq2}
\end{equation}
\end{theorem}

\section{Some examples}

\begin{example}
\rm
Let $G(z,w) = \alpha(w) z + \beta(w)$ with $\beta(0)=0$.
Then \reff{eq.star3} tells us that
\begin{equation}
   \varphi(w)
   \;=\;
   \sum_{m=1}^\infty {1 \over m} \,
        [\zeta^{m-1}] \, [\alpha(w) \zeta + \beta(w)]^m
   \;=\;
   \sum_{m=1}^\infty \alpha(w)^{m-1} \beta(w)
   \;,
\end{equation}
which sums to the correct answer $\beta(w) / [1 -\alpha(w)]$
provided that $|\alpha(w)| < 1$.
So some condition like $|(\partial G/\partial z)(0,0)| < 1$
is needed in order to ensure convergence of the series \reff{eq.star3}.
\qed
\end{example}

\begin{example}
 \label{exam3.2}
\rm
The inverse-function problem $f(z)=w$
with $f(z) = \sum_{n=1}^\infty a_n z^n$ ($a_1 \neq 0$)
can be written in the form $z = G(z,w)$ in a variety of different ways.
The most obvious choice is
\begin{equation}
   G(z,w)  \;=\;  {z \over f(z)} \, w \;,
\end{equation}
which leads to the usual form of the Lagrange inversion formula:
\begin{equation}
   f^{-1}(w)
   \;=\;
   \sum_{m=1}^\infty {w^m \over m} \, a_1^{-m} \,
            [\zeta^{m-1}]
            \biggl( 1 + \sum\limits_{n=2}^\infty {a_n \over a_1} \zeta^{n-1}
            \biggr)^{\! -m}
   \;.
 \label{eq.lagrange.usual}
\end{equation}
An alternative choice, proposed by Yuzhakov \cite{Yuzhakov_75a}, is
\begin{equation}
   G(z,w)  \;=\;  {w \over a_1} \,-\, {f(z) - a_1 z  \over a_1} \;,
 \label{eq.Gzw.yuzhakov}
\end{equation}
which leads to
\begin{eqnarray}
   f^{-1}(w)
   & = &
   \sum_{m=1}^\infty {a_1^{-m} \over m}
   \sum_{\ell=0}^{m} \binom{m}{\ell} w^{m-\ell} (-1)^\ell \,
            [\zeta^{m-1}]
            \biggl( \sum\limits_{n=2}^\infty a_n \zeta^n \biggr)^{\! \ell}
   \;.
 \label{eq.lagrange.yuzhakov}
\end{eqnarray}
After some straightforward algebra involving binomial coefficients,
both forms can be shown to yield the same result:
\begin{equation}
   f^{-1}(w)
   \;=\;
   \sum_{m=1}^\infty {w^m \over m!}
   \!\!
   \sum_{\substack{
            k_2,k_3,\ldots \ge 0 \\[1mm]
            \sum (n-1)k_n = m-1
        }}
   \!\!
   (-1)^{\sum k_n} \: a_1^{-(1 + \sum nk_n)} \:
   {\textstyle \bigl( \sum nk_n \bigr)!
    \over
    \, k_2! \, k_3! \, \cdots \,
   }
   \:
   \prod_{n=2}^\infty a_n^{k_n}
 \label{eq.inverse.explicit}
\end{equation}
(see also \cite{Gessel_88}).
However, \reff{eq.lagrange.usual}/\reff{eq.inverse.explicit}
expands $f^{-1}(w)$ as a power series in $w$,
while \reff{eq.lagrange.yuzhakov} expands $f^{-1}(w)$ as a series
in a different set of polynomials in $w$.
%
%

Consider, for instance, $f(z) = z e^{-z}$.
Then the usual Lagrange inversion formula
\reff{eq.star1}
gives
\begin{equation}
   f^{-1}(w)
   \;=\;
   \sum_{m=1}^\infty m^{m-1} \, {w^m \over m!}
   \;,
\end{equation}
which is well known \cite[Propositions~5.3.1 and 5.3.2]{Stanley_99}
to be the exponential generating function for rooted trees
(i.e., there are $m^{m-1}$ distinct rooted trees on $m$ labeled vertices).
On the other hand, Yuzhakov's version
\reff{eq.Gzw.yuzhakov}/\reff{eq.lagrange.yuzhakov} gives,
after a short calculation, the alternate representation
\begin{equation}
   f^{-1}(w)
   \;=\;
   \sum_{m=1}^\infty  P_m(w)
 \label{eq.tree.yuzhakov.1}
\end{equation}
where
\begin{equation}
   P_m(w)  \;=\;  (-1)^{m-1} \sum_{k = \lceil (m+1)/2 \rceil}^m
       {(m-1)! \over k! \, (k-1)!} \, \stirlingsubset{k-1}{m-k} \, w^k
   \;;
 \label{eq.tree.yuzhakov.2}
\end{equation}
here $\stirlingsubset{n}{k}$ denotes the Stirling subset numbers
(also known as Stirling numbers of the second kind),
i.e.\ the number of partitions of an $n$-element set
into $k$ nonempty blocks \cite{Graham_94}.\footnote{
   The key step in the derivation of
   \reff{eq.tree.yuzhakov.1}/\reff{eq.tree.yuzhakov.2}
   is the well-known identity \cite[eq.~(7.49)]{Graham_94}
   $$ (e^{-z}-1)^k  \;=\;
      k! \sum_{n=k}^\infty \stirlingsubset{n}{k} \, {(-z)^n \over n!}  \;. $$
}
I wonder whether the coefficients in $P_m(w)$ have any combinatorial meaning.
\qed
\end{example}

\begin{example}
 \label{exam3.3}
\rm
Here is an application from my own current research \cite{Sokal_Fxy}.
The function
\begin{equation}
   F(x,w) \;=\; \sum\limits_{n=0}^\infty {x^n \over n!} \, w^{n(n-1)/2}
 \label{def.F}
\end{equation}
arises in enumerative combinatorics in the generating function
for the Tutte polynomials of the complete graphs $K_n$
\cite{Tutte_67,Scott-Sokal_expidentities}
and in statistical mechanics as the grand partition function
of a single-site lattice gas with fugacity $x$ and
two-particle Boltzmann weight $w$ \cite{Scott-Sokal_lovasz}.
Let us consider $F$ as a function of complex variables $x$ and $w$
satisfying $|w| \le 1$:
it is jointly analytic on $\C \times \D$
and jointly continuous on $\C \times \Dbar$
(here $\D$ and $\Dbar$ denote the open and closed unit discs in $\C$,
respectively),
and it is an entire function of $x$ for every $w \in \Dbar$.

When $w=0$, we have $F(x,0) = 1+x$, which has a simple zero at $x=-1$.
One therefore expects --- and can easily prove using Rouch\'e's theorem ---
that for small $|w|$ there is a unique root of $F(x,w)$ near $x=-1$,
which can be expanded in a convergent power series
\begin{equation}
   x_0(w)  \;=\;  -1 \,-\, \sum_{n=1}^\infty a_n w^n
   \;.
 \label{def.series.x0}
\end{equation}
The coefficients $\{a_n\}$ can of course be computed by substituting the
series \reff{def.series.x0} into \reff{def.F} and equating
term-by-term to zero;  but a more efficient method is to use
the implicit function formula \reff{eq.star3}.
It suffices to set $x=-1-z$ and define
\begin{equation}
   G(z,w) \;=\; \sum\limits_{n=2}^\infty {(-1-z)^n \over n!} \, w^{n(n-1)/2}
   \;.
\end{equation}
We then obtain
\begin{eqnarray}
   -x_0(w)
   & = &
   1 \,+\, \smfrac{1}{2} w \,+\, \smfrac{1}{2} w^2
     \,+\, \smfrac{11}{24} w^3 \,+\, \smfrac{11}{24} w^4
     \,+\, \smfrac{7}{16} w^5 \,+\, \smfrac{7}{16} w^6
     \,+\, \smfrac{493}{1152} w^7 \,+\, \smfrac{163}{384} w^8
   \nonumber \\[1mm]
   &  &  \quad
     \,+\, \smfrac{323}{768} w^9 \,+\, \smfrac{1603}{3840} w^{10}
     \,+\, \smfrac{57283}{138240} w^{11} \,+\, \smfrac{170921}{414720} w^{12}
     \,+\, \ldots
   \;.
 \label{eq.series.x0}
\end{eqnarray}
I conjecture --- but have thus far been unable to prove ---
that all the coefficients in this power series are nonnegative.
Since it is known \cite{Morris_72,Liu_98,Langley_00}
that $x_0(w)$ is analytic in a complex neighborhood of the
real interval $0 < w < 1$,
this conjecture would imply,
by the Vivanti--Pringsheim theorem \cite[Theorem 5.7.1]{Hille_73},
that $x_0(w)$ is in fact analytic in the whole disc $|w| < 1$,
i.e.\ that the series \reff{eq.series.x0}
has radius of convergence exactly 1.
(It is not hard to show that $x_0(w) \to -\infty$ as $w \uparrow 1$,
 so that the radius of convergence cannot be bigger than 1.)
   
Using \reff{eq.star4} we can also compute power series for functions
of $x_0(w)$.  For instance, we have
\begin{eqnarray}
   \log[-x_0(w)]
   & = &
           \smfrac{1}{2} w \,+\, \smfrac{3}{8} w^2
     \,+\, \smfrac{1}{4} w^3 \,+\, \smfrac{41}{192} w^4
     \,+\, \smfrac{13}{80} w^5 \,+\, \smfrac{85}{576} w^6
     \,+\, \smfrac{83}{672} w^7 \,+\, \smfrac{227}{2048} w^8
   \nonumber \\[1mm]
   &  &  \quad
     \,+\, \smfrac{2065}{20736} w^9 \,+\, \smfrac{4157}{46080} w^{10}
     \,+\, \smfrac{6953}{84480} w^{11} \,+\, \smfrac{252449}{3317760} w^{12}
     \,+\, \ldots
   \;.
 \label{eq.series.logx0}
\end{eqnarray}
I conjecture that all the coefficients in \reff{eq.series.logx0}
are nonnegative.
By exponentiation this implies the preceding conjecture, but is stronger.
We also have
\begin{eqnarray}
   - \, {1 \over x_0(w)}
   & = &
   1 \,-\, \smfrac{1}{2} w \,-\, \smfrac{1}{4} w^2
     \,-\, \smfrac{1}{12} w^3 \,-\, \smfrac{1}{16} w^4
     \,-\, \smfrac{1}{48} w^5 \,-\, \smfrac{7}{288} w^6
     \,-\, \smfrac{1}{96} w^7 \,-\, \smfrac{7}{768} w^8
   \nonumber \\
   &  &  \quad
     \,-\, \smfrac{49}{6912} w^9 \,-\, \smfrac{113}{23040} w^{10}
     \,-\, \smfrac{17}{4608} w^{11} \,-\, \smfrac{293}{92160} w^{12}
     \,-\, \ldots
   \;.
 \label{eq.series.1overx0}
\end{eqnarray}
I conjecture that all the coefficients in \reff{eq.series.1overx0}
after the constant term are nonpositive.
This implies the preceding two conjectures, but is even stronger.
Using {\sc Mathematica} I have verified all three conjectures
through order $w^{60}$.
Indeed, by exploiting the connection between $F(x,w)$
and the generating polynomials $C_n(v)$ of connected graphs
on $n$ labeled vertices \cite{Tutte_67,Scott-Sokal_expidentities,Sokal_Fxy}
--- or equivalently the inversion enumerator for trees, $I_n(w)$
\cite[Exercise~5.48, pp.~93--94 and 139--140]{Stanley_99} ---
I have computed the series $x_0(w)$ and verified these conjectures
through order $w^{775}$.
%
%

The relative simplicity of the coefficients in \reff{eq.series.1overx0}
compared to \reff{eq.series.x0}/\reff{eq.series.logx0}
suggests that $-1/x_0(w)$ may have a simpler combinatorial interpretation
than $-x_0(w)$ or $\log[-x_0(w)]$.
Please note also that since $x_0(w) \to -\infty$ as $w \uparrow 1$,
the coefficients $\{b_n\}$ in
$-1/x_0(w) = 1 - \sum_{n=1}^\infty b_n w^n$
--- if indeed they are nonnegative --- add up to~1,
so they are the probabilities for a positive-integer-valued random variable.
What might such a random variable be?
Could this approach be used to prove the nonnegativity of $\{b_n\}$?

Please note also that the coefficients $\{a_n\}$ and $\{b_n\}$
satisfy the discrete-time renewal equation
\be
   a_n  \;=\;  \sum_{k=1}^{n-1} b_k a_{n-k}
   \;.
\ee
Therefore, if the $\{b_n\}$ are nonnegative,
they can be interpreted \cite[Chapter~XIII]{Feller_68}
as the probability distribution for first occurrences
(or equivalently for waiting times between successive occurrences)
of a recurrent event $\scre$, while the $\{a_n\}$ are the probabilities
of occurrence {\em tout court}\/:
\begin{subeqnarray}
   a_n  & = &  {\mathbb P}(\hbox{$\scre$ occurs at the $n$th trial}) \\
   b_n  & = &  {\mathbb P}(\hbox{$\scre$ occurs for the first time
                                                at the $n$th trial}) \qquad
\end{subeqnarray}
What might such a family of recurrent events be?

For what it's worth, if we define $c_N = 1 - \sum_{n=1}^N b_n$,
we find empirically (at least up to $N=775$) that
%
%
that $c'_N = 2 \, N! \, c_N$ is an integer, with
\begin{eqnarray}
   & & c'_1, \ldots, c'_{20} \;=\;
   1,\, 1,\, 2,\, 5,\, 20,\, 85,\, 490,\, 3185,\, 23520,\, 199605,\, 1901130,\,
         \nonumber \\
   & & \quad
   19767825,\, 223783560,\, 2806408605,\, 37447860450,\, 540137222625,\,
         \nonumber \\
   & & \quad
   8284392916800,\, 135996789453525,\, 2363554355812650,\, 43437044503677825
         \nonumber \\  \\[-12mm] \nonumber
\end{eqnarray}
Can anyone figure out a combinatorial interpretation of these numbers?

It is also known \cite{Sokal_Fxy} that $-x_0(w)$ diverges as $w \uparrow 1$
with leading term $e^{-1} (1-w)^{-1}$,
which suggests that we have
$\lim_{n\to\infty} a_n = e^{-1}$
and $\sum_{n=1}^\infty nb_n = e$.
If we define $d_N = \sum_{n=1}^N nb_n$
and $e_N = \sum_{n=1}^N [1/(n-1)! \,-\, nb_n]$,
we find empirically (at least up to $N=775$) that
$d'_N = 2 \, (N-1)! \, d_N$ and $e'_N = 2 \, (N-1)! \, e_N$ are integers, with
\begin{eqnarray}
   & & d'_1, \ldots, d'_{20} \;=\;
   1,\, 2,\, 5,\, 18,\, 77,\, 420,\, 2625,\, 19110,\, 158025,\, 1457820,\,
      14872725,\,
         \nonumber \\
   & & \quad
   166645710,\, 2032946685,\, 26754868140,\, 379216422585,\, 5747274883350,\,
         \nonumber \\
   & & \quad
   92854338001425,\, 1591646029073100,\, 28870013167120125,\,
      552364292787857550
         \nonumber \\ \\
   & & e'_1, \ldots, e'_{20} \;=\;
   1,\, 2,\, 5,\, 14,\, 53,\, 232,\, 1289,\, 8290,\, 61177,\, 515000,\,
      4855477,\,
         \nonumber \\
   & & \quad
   50364514,\, 571176005,\, 7098726832,\, 94733907025,\, 1361980060802,\,
         \nonumber \\
   & & \quad
   20893741105009, 342071315736280, 5936899039448717, 108967039136950450
         \nonumber \\  \\[-12mm] \nonumber
\end{eqnarray}
Might these numbers have some combinatorial interpretation?
\qed
\end{example}

\section{Possible multidimensional extensions}

In this paper I have restricted attention to the implicit-function
problem in one complex variable
(i.e., $z \in \C$ though $w$ lies in an arbitrary space $W$).
Yuzhakov \cite{Yuzhakov_75a,Yuzhakov_75b,Aizenberg_83}
goes much farther:  he gives a beautiful explicit formula for the
solution of the multidimensional implicit-function problem $F(z,w) = 0$
with $z \in \C^N$, $w \in \C^M$ and $F \colon\, \C^N \times \C^M \to \C^N$
under the usual hypothesis that
the linear operator $(\partial F/\partial z)(0,0)$ is nonsingular.
I suspect that the approach of the present paper can likewise
be extended to the corresponding multidimensional situation ---
that is, $z = G(z,w)$ with $z \in \C^N$, $w \in W$ and
$G \colon\, \C^N \times W \to \C^N$ ---
under the hypothesis that the linear operator $(\partial G/\partial z)(0,w_0)$
has spectral radius $< 1$.
Indeed, such a result presumably holds when $\C^N$
is replaced by a complex Banach space.
As we have seen, the proof of Theorem~\ref{thm1} given here
applies verbatim when $w$ lies in an arbitrary space,
since the variables $w$ simply ``go for the ride''.
But multidimensional $z$ is a genuine generalization;
and for lack of time and competence, I have not attempted to pursue it.
See \cite{Good_60,Good_65,Henrici_84,Labelle_85a,Labelle_85b,Gessel_87,%
Bergeron_98,Abdesselam_03,Bousquet_03}
for information on multidimensional Lagrange inversion formulae,
and \cite{Krantz_02} for a survey of implicit function theorems.

\section*{Acknowledgments}

I wish to thank Alex Eremenko, Ira Gessel, Steven Krantz,
Gilbert Labelle, Pierre Leroux and Richard Stanley
for many helpful comments.
In particular, Pierre Leroux gave comments on an early draft of this paper
and kindly drew my attention to work on
Lagrange inversion and implicit-function formulae
from the point of view of the theory of combinatorial species
\cite{Joyal_81,Labelle_85a,Labelle_85b,Gessel_95,Bergeron_98}.
I dedicate this paper to his memory.

I also wish to thank the Institut Henri Poincar\'e -- Centre Emile Borel
for hospitality during the programme on
Interacting Particle Systems, Statistical Mechanics and Probability Theory
(September--December 2008),
where this work was (almost) completed.

This research was supported in part
by U.S.\ National Science Foundation grant PHY--0424082.


\end{document}